\documentclass[11pt, psamsfonts]{amsart}
\usepackage[foot]{amsaddr}
\usepackage{amssymb,amsfonts}
\usepackage[all,arc]{xy}
\usepackage{enumerate}
\usepackage{mathrsfs}
\usepackage{tikz-cd, tikz}
\usepackage[hidelinks]{hyperref}

\usepackage[
backend=biber,
style=numeric,
]{biblatex}

\newtheorem{thm}{Theorem}[section]
\newtheorem{cor}[thm]{Corollary}
\newtheorem{prop}[thm]{Proposition}
\newtheorem{lem}[thm]{Lemma}

\theoremstyle{definition}
\newtheorem{defn}[thm]{Definition}

\theoremstyle{remark}
\newtheorem{rem}[thm]{Remark}

\makeatletter
\let\c@equation\c@thm
\makeatother
\numberwithin{equation}{section}

\title{Further Notes On Tightness}

\author[Khlaif]{Nasief Khlaif}
\author[Saleh]{Mohammad Saleh}
\email{nkhlaif@birzeit.edu}
\email{msaleh@birzeit.edu}
\address{Birzeit University, Birzeit, Palestine}

\begin{document}

\begin{abstract}
    Tight and essentially tight modules generalize weakly injective modules. Essential tightness requires embeddings to be essential. This restriction makes the two notions fundamentally different. In this note, we investigate cases when these two notions coincide. Moreover, we examine the cases when essentiality is imposed only on one of the embeddings rather than both. This allows us to define a special class of tight and essentially tight modules and a generalization of both.
\end{abstract}
\maketitle

\section{Introduction}
Weakly injective modules were introduced by Jain and Lopez \cite{jain1990rings} as a generalization of injective modules in their study of CEP rings, which are rings whose cyclics are essentially embeddable in a projective module. They proved that a ring $R$ is CEP if and only if it is artinian and every indecomposable projective module is uniform and weakly $R-$injective. Later, Golan and Lopez \cite{golan1991qi} further generalized weak injectivity into what they coined as tightness. Saleh \cite{saleh1999note} proved a more general result than Jain and Lopez's characterization of CEP rings replacing weak injectivity by tightness. Recently, Koşan, Quynh and Serap \cite{kocsan2017weakly} introduced another generalization to weak injectivity, namely essential tightness. They proved a characterization of CEP ring similar to Saleh's using essential tightness.

It is noted in \cite{kocsan2017weakly} that not every tight module is essentially tight. In this short note, we provide different cases when these notions are actually equivalent. To do this, we review, in Section 2, fundamental results on weak injectivity, tightness and essential tightness from the literature, covering fundamentals of \textit{q.f.d.} and semiprime Goldie rings. The equivalence between tightness and essential tightness is discussed in Section 3. The main theorem in this regard is Theorem 3.11. 

Both tightness and essentially tight modules generalize weak injectivity by imposing or removing essentiality on both embeddings. Imposing it on only one of the embeddings provides us with two further generalizations of weak injectivity: one is a special case of both tight and essentially tight, and the other is a generalization of both. In Section 4, we define these and prove some results concerning them, focusing on CEP ring characterizations. It turns out that most results are similar to those for weakly injective, tight and essentially tight modules. Throughout this paper, it is assumed that all rings are associative with unity and all modules are right and unital.
\section{Preliminaries}
A module $M$ is called weakly $N-$injective if for any $\sigma: N \to E(M)$, there exists a submodule $X$ of $E(M)$ such that $\sigma(N) \subset X \simeq M$. $M$ is called weakly-injective if it is weakly $N-$injective for every finitely generated module $N$. On the other hand, $M$ is said to be (essentially) $N-$tight if for every submodule $K$ of $N$, whenever $N/K$ (essentially) embeds into $E(M)$, then $N/K$ (essentially) embeds into $M$. When $M$ is (essentially) $N-$tight for every finitely generated module $N$, then $M$ is said to be (essentially) tight. Hence, any weakly injective module is both tight and essentially tight. The following characterization appears in \cite[Proposition 4.1]{jain1990rings} and \cite[Theorem 3.2]{kocsan2017weakly}, which guarantees the existence of monomorphisms from quotients by submodules.

\begin{prop}
\label{21}
    Given $R-$modules $M$ and $N$, the following are equivalent.
    \begin{enumerate}
        \item $M$ is weakly $N-$injective (essentially $N-$tight).
        \item $M$ is weakly $N/K-$injective (essentially $N/K-$tight) for any $K \subset N$.
        \item For every submodule $K$ of $N$ and for every (essential) monomorphism $h:N/K \to E(M)$ there exist (essential) monomorphisms $\sigma:M \to E(M)$ and $g: N/K \to M$ such that $h = \sigma \circ g$.
    \end{enumerate}
\end{prop}

One important result that relates tight and weakly injective modules is the following result by Jain and Lopez \cite[Proposition 2.2]{jain1995weakly}. This result along with results presented next allow characterizing different cases when tightness and weak injectivity are the same.

\begin{prop}
    If $M$ is a module with $E(M)$ is a direct sum of indecomposables. Then $M$ is tight if and only if $M$ is weakly injective.
\end{prop}

\begin{rem}
    For any module $M$ over a noetherian ring, $E(M)$ can be written as a direct sum of indecomposables.
\end{rem}

A module is said to be uniform if any nonzero submodule is essential. A module $M$ is said to have a Goldie (uniform) dimension $u.dim(M) = n \in \mathbf{Z}$ if there exists a finite set of uniform submodules $U_i$ such that $\oplus_{i=1}^{n}U_{i}$ is an essential submodule of $M$. A ring $R$ is said to be \textit{q.f.d.} if every cyclic module has a finitely generated, possibly trivial, socle. Hence, a ring is \textit{q.f.d.} if every cyclic module has a finite Goldie dimension. One important class of \textit{q.f.d.} rings is the class of noetherian rings. One important result over \textit{q.f.d.} rings is that of Lopez \cite[Theorem 3.1]{lopez1992rings}, which establishes an equivalence between weak injectivity and tightness over \textit{q.f.d.} rings.

\begin{thm}
    Over a q.f.d. ring, every tight module is weakly injective.
\end{thm}

The following result relates weak injectivity, tightness and essential tightness over \textit{q.f.d.} rings from both \cite[Theorem 2.6]{jain2018survey} and \cite[Theorem 3.10]{kocsan2017weakly}.
\begin{thm}
    The following are equivalent for a ring $R$.
    \begin{enumerate}
        \item $R$ is \textit{q.f.d.}
        \item Every direct sum of injectives is weakly-injective.
        \item Every direct sum of injectives is tight.
        \item Every direct sum of tight is tight.
        \item Every direct sum of weakly-injectives is tight.
        \item Every direct sum of weakly-injectives is $R-$tight.
        \item Every direct sum of indecomposable injectives is $R-$tight.
        \item Every direct sum of essentially tight is essentially tight.
        \item Every direct sum of injective $R-$modules is essentially tight.
        \item Every direct sum of indecomposable injective $R-$modules is essentially tight.
    \end{enumerate}
\end{thm}

A right \textit{Goldie ring} is a ring with a finite Goldie dimension as a module over itself, $u.dim \,\,R_R < \infty$, and every ascending chain of annihilator ideals stabilizes. Similarly, left Goldie rings can be defined. A ring is Goldie if it is both left and right Goldie. One important result that plays a pivotal role in this note is the following result from \cite[Proposition 3.8]{lopez1993some}. It provides a case when specific modules are always weakly injective.

\begin{prop}
    Every nonsingular module over a semiprime Goldie ring $R$ is weakly injective.
\end{prop}

\section{Tight Meets Essentially Tight}
Tightness and essential tightness are two generalizations of weak injectivity. Although they satisfy identical results, they are two different classes of modules. Below, we answer a question raised by Koşan \cite{kocsan2017weakly} concerning equivalence between tight and essentially tight modules. Motivated by results in the previous section, the following results provide connections between tightness and essential tightness and their proofs are similar.

\subsection{Specific Types of Modules}
Some specific types of modules lie at the intersection of tight and essentially tight modules. Some of those are presented here as a further extension of results relating tightness and weak injectivity.

\begin{prop}
    If $M$ is a module with $E(M)$ is a direct sum of indecomposables. Then $M$ is essentially tight if and only if $M$ is weakly injective.
\end{prop}

\begin{cor}
    Over a noetherian ring, a module is weakly injective if and only if it is essentially tight.
\end{cor}
\begin{proof}
    This follows from Remark 2.3.
\end{proof}

\begin{cor}
    If $M$ is a module with $E(M)$ is a direct sum of indecomposables. Then, the following are equivalent:
    \begin{enumerate}
        \itemsep0em
        \item $M$ is tight.
        \item $M$ is weakly injective.
        \item $M$ is essentially tight.\hfill \qed
    \end{enumerate}
\end{cor}

\begin{cor}
    Over a noetherian ring, the following are equivalent:
    \begin{enumerate}
        \itemsep0em
        \item $M$ is tight.
        \item $M$ is weakly injective.
        \item $M$ is essentially tight.\hfill \qed
    \end{enumerate}
\end{cor}

\begin{prop}
    Let $M$ be a uniform module. $M$ is $N-$tight if and only if $M$ is essentially $N-$tight.
\end{prop}
\begin{proof}
    Let $M$ be a uniform module and $N$ be any other module with $K \subset N$. Assume $M$ is essentially $N-$tight. Let $f: N/K \to E(M)$ be an embedding. Then, since $M$ is uniform, $E(M)$ is also uniform and hence $f$ is an essential embedding. Since $M$ is essentially $N-$tight, we have an essential embedding $\tilde{f}: N/K \to M$, hence $M$ is $N-$tight.
    
    Conversely, let $M$ be $N-$tight. Then for the embedding $g: N/K \to E(M)$, it is necessarily essential since $E(M)$ is uniform and any nonzero submodule is essential, thus $g(N/K) \trianglelefteq E(M)$. By the $N-$tightness of $M$, we have an embedding $\tilde{g}: N/K \to M$, which is also essential by the uniformity of $M$. Hence $M$ is essentially $N-$tight.
\end{proof}
\begin{cor}
    For a uniform module, the following are equivalent:
    \begin{enumerate}
        \itemsep0em
        \item $M$ is tight.
        \item $M$ is weakly injective.
        \item $M$ is essentially tight.\hfill \qed
    \end{enumerate}
\end{cor}

\subsection{Over Some Classes of Rings}
Some classes of rings provide rich module structure that was investigated with respect to tight and weakly injective modules as stated in Section 2. Most of the results apply naturally in a similar manner without much modification.

\begin{thm}
    Essentially tight modules are weakly injective over \textit{q.f.d.}'s.
\end{thm}

\begin{thm}
    Over a \textit{q.f.d.} ring, the following are equivalent: 
    \begin{enumerate}
        \itemsep0em
        \item $M$ is tight.
        \item $M$ is weakly injective.
        \item $M$ is essentially tight.\hfill \qed
    \end{enumerate}
\end{thm}

\begin{cor}
    Over a semiprime Goldie ring, every nonsingular essentially tight module is weakly injective.\hfill\qed
\end{cor}

\begin{cor}
    For nonsingular modules over semiprime Goldie rings, the following are equivalent:
    \begin{enumerate}
        \itemsep0em
        \item $M$ is tight.
        \item $M$ is weakly injective.
        \item $M$ is essentially tight.\hfill \qed
    \end{enumerate}
\end{cor}

This section concludes with the following theorem, which is one of the main results of this note, summing up all results presented thus far on the equivalence between tight and essentially tight modules.
\begin{thm}
    Essential tightness of a module $M$ is equivalent to tightness whenever one of the following cases holds.
    \begin{enumerate}
        \itemsep0em
        \item $M$ is uniform.
        \item $E(M)$ can be written as a direct sum of indecomposables.
        \item Over a \textit{q.f.d.} ring.
        \item $M$ is nonsingular over a semiprime Goldie ring.
    \end{enumerate}
\end{thm}

\section{New Classes}

In defining essential tightness, essentiality is imposed on both embeddings of tightness. Imposing essentiality on one of the embeddings allows to define two new types of modules. The first is a special case of both tight and essentially tight modules, coined as minimally tight. On the other hand, the second is a generalization of both, coined as quasi tight. In this section, we study basic properties of those types and prove some results concerning them.

\subsection{Quasi Tight Modules}
Quasi tight modules generalize both tight and essentially tight modules. They provide a more general characterization of CEP rings in spirit of the characterizations  of CEP rings with respect to weakly injective, tight and essentially tight modules.

\begin{defn}[Relatively Quasi Tight Module]
A module $M$ is quasi $N-$tight if for any submodule $K$ of $N$, $N/K$ essentially embeds in $E(M)$ implies that $N/K$ embeds in $M$.    
\end{defn}

\begin{defn}[Quasi Tight Module]
    A module $M$ is said to be quasi tight if it is quasi $N-$tight for every finitely generated module $N$.
\end{defn}

\begin{cor}
    Every tight module is quasi tight.
\end{cor}
\begin{proof}
    Let $M$ be an $N-$tight module, assume $\phi: N/K \to E(M)$ be an essential embedding for a submodule $K$ of $N$. Hence $\phi$ is an embedding of $N/K$ in $E(M)$. Since, $M$ is $N-$tight, then $N/K$ is embeddable into $M$. Hence $M$ is quasi $N-$tight.
\end{proof}

\begin{cor}
    Every essentially tight module is quasi tight.
\end{cor}
\begin{proof}
    Let $M$ be an essentially $N-$tight module. Assume $\phi: N/K \to E(M)$ is an essential embedding for a submodule $K$ of $N$. Since $M$ is essentially $N-$tight, then $N/K$ is essentially embeddable into $M$ and so embeddable into $M$. So, $M$ is quasi $N-$tight.
\end{proof}

\begin{cor}
    Any weakly injective module is quasi tight \hfill\qed
\end{cor}

\begin{rem}
    Let $M$ be a finitely generated module. $M$ is weakly injective if and only if it is quasi tight.
\end{rem}

Quasi tight modules turn out to satisfy some properties similar to those of (essentially) tight modules. The following results provide some insights into the behavior of quasi tight modules in different settings.

\begin{lem}
    For modules $M$ and $N$. The following are equivalent:
    \begin{enumerate}
        \item $M$ is quasi $N-$tight.
        \item $M$ is quasi $N/K-$tight for any $K \subset N$.
        \item For every submodule $K$ of $N$ and for every essential monomorphism $h:N/K \to E(M)$ there exist an essential monomorphism $\sigma:M \to E(M)$ and a monomorphism $g: N/K \to M$ such that $h = \sigma \circ g$.
    \end{enumerate}
\end{lem}

\begin{proof}
    Similar to proof of Theorem 3.2 in Koşan \cite{kocsan2017weakly}.
\end{proof}

\begin{thm}
    Let $M$ be quasi $N-$tight. Then
    \begin{enumerate}
            \item If $K \trianglelefteq N$, then $M$ is quasi $K-$tight.
        \item If $K \cong N$, then $M$ is quasi $K-$tight.
        \item If $L$ is an essential extension of $M$, then $L$ is quasi $N-$tight.
    \end{enumerate}
\end{thm}

\begin{proof}
    Similar to proof of Theorem 2.5 in Koşan \cite{kocsan2017weakly}.
\end{proof}

\begin{thm}
    The following are true over any ring $R$.
    \begin{enumerate}
        \item Let $N$ be an uniform module. If $P$ and $Q$ are quasi $N-$tight, then so is $P \oplus Q$.
        \item Quasi tight modules are closed under finite direct sums.
    \end{enumerate}
\end{thm}

\begin{proof}
    Similar to proof of Theorem 3.5 in Koşan \cite{kocsan2017weakly}.
\end{proof}

Since the main constructions of the different types of modules so far presented originated from characterizing CEP rings, it would be worthwhile to characterize CEP rings in terms of quasi tightness. This is the focus of the next results. For more compact notation, two modules $M$ and $N$ are said to be \textit{relatively quasi tight} if $M$ is quasi $N-$tight and $N$ is quasi $M-$tight.

\begin{lem}
    Let $N$ and $M$ be finitely generated modules over an artinian ring $R$. If $M$ and $N$ are relatively quasi tight and $\mathbf{Soc}(M) \cong \mathbf{Soc}(N)$, then $M \cong N$.
\end{lem}

\begin{proof}
    Similar to proof of Lemma 1 in Saleh \cite{saleh1999note}.
\end{proof}

\begin{lem}
    Let $R$ be an artinian ring such that all indecomposable projective $R-$modules are uniform and quasi $R-$tight. Then
    \begin{enumerate}
        \item Every simple module is isomorphic to the socle of an indecomposable projective module,
        \item Every simple $R-$module is embeddable in $\mathbf{Soc}(R)$, and
        \item If $P$ and $Q$ are projectives with $\mathbf{Soc}(P) \cong \mathbf{Soc}(Q)$ then $P \cong Q$.
    \end{enumerate}
\end{lem}

\begin{proof}
    Similar to proof of Lemma 3 in Saleh \cite{saleh1999note}.
\end{proof}

\begin{thm}
    A ring $R$ is CEP if and only if 
        \begin{enumerate}
    \itemsep0em
        \item $R$ is artinian.
        \item Indecomposable projective right $R-$modules are uniform and quasi $R-$tight.
    \end{enumerate}
\end{thm}

\begin{proof}
    Similar to proof of Theorem 1 in Saleh \cite{saleh1999note} and Theorem 5 in \cite{jain1990rings}.
\end{proof}

In the spirit of the study presented on \textit{q.f.d.} rings, results extend from weakly injective and (essentially) tight modules to quasi tight modules naturally along the same lines. This allows us to further extend the equivalence results for \textit{q.f.d.} rings.

\begin{thm}
The following are equivalent for a ring $R$.
    \begin{enumerate}
        \item $R$ is \textit{q.f.d.}.
        \item Every direct sum of quasi tight modules is quasi tight.
        \item Every direct sum of injective $R-$modules is quasi tight.
        \item Every direct sum of indecomposable injectives is quasi tight.
\end{enumerate}
\end{thm}

\begin{proof}
    Similar to Theorem 3.10 in Koşan \cite{kocsan2017weakly} and Theorem 2.6 in \cite{jain2018survey}.
\end{proof}

\subsection{Minimally Tight Modules}
Such modules represent a special class of both tight and essentially tight modules at the border line between weakly injective and (essentially) tight modules. They share many results with those classes. In this section, we define those modules and present some of their properties. Proofs of different results presented here are omitted as they go along the same lines with previous results.

\begin{defn}[Relatively Minimally Tight Module]
A module $M$ is minimally $N-$tight if for any submodule $K$ of $N$, $N/K$ embeds in $E(M)$ implies that $N/K$ essentially embeds in $M$.    
\end{defn}

\begin{defn}[Minimally Tight Module]
    A module $M$ is said to be minimally tight if it is minimally $N-$tight for every finitely generated module $N$.
\end{defn}

\begin{cor}
    Every minimally tight module is tight.
\end{cor}
\begin{proof}
    Let $M$ be a minimally $N-$tight module. Assume $\phi: N/K \to E(M)$ is an embedding for a submodule $K$ of $N$. Since $M$ is minimally $N-$tight, then $N/K$ is essentially embeddable into $M$, hence embeddable into $M$. Therefore, $M$ is $N-$tight.
\end{proof}

\begin{cor}
    Every minimally tight module is essentially tight.
\end{cor}
\begin{proof}
    Let $M$ be a minimally $N-$tight module. Assume $\phi: N/K \to E(M)$ is an embedding for a submodule $K$ of $N$. Since $M$ is minimally $N-$tight, then $N/K$ is essentially embeddable into $M$. So $M$ is essentially $N-$tight.
\end{proof}

\begin{lem}
    For modules $M$ and $N$. The following are equivalent:
    \begin{enumerate}
        \item $M$ is minimally $N-$tight.
        \item $M$ is minimally $N/K-$tight for any $K \subset N$.
        \item For every submodule $K$ of $N$ and for every monomorphism $h:N/K \to E(M)$ there exist a monomorphism $\sigma:M \to E(M)$ and an essential monomorphism $g: N/K \to M$ such that $h = \sigma \circ g$. \hfill\qed
    \end{enumerate}
\end{lem}

\begin{cor}
    Every weakly injective module is minimally tight.\hfill\qed
\end{cor}

Similar to quasi tight modules characterizations, minimally tight modules also provide a characterization of CEP rings. Two modules $M, N$ are said to be \textit{relatively minimally tight} if $M$ is minimally $N-$tight and $N$ is minimally $M-$tight. Proofs also go along the same lines, hence they are omitted.

\begin{lem}
    Let $N$ and $M$ be finitely generated modules over an artinian ring $R$. If $M$ and $N$ are relatively minimally tight and $\mathbf{Soc}(M) \cong \mathbf{Soc}(N)$, then $M \cong N$.\qed
\end{lem}

\begin{lem}
    Let $R$ be an artinian ring such that all indecomposable projective $R-$modules are uniform and minimally $R-$tight. Then
    \begin{enumerate}
        \item Every simple module is isomorphic to the socle of an indecomposable projective module,
        \item Every simple $R-$module is embeddable in $\mathbf{Soc}(R)$, and
        \item If $P$ and $Q$ are projective with $\mathbf{Soc}(P) \cong \mathbf{Soc}(Q)$, then $P \cong Q$.\qed
    \end{enumerate}
\end{lem}

\begin{thm}
    A ring $R$ is CEP if and only if 
        \begin{enumerate}
    \itemsep0em
        \item $R$ is artinian.
        \item Indecomposable projectives are uniform and minimally $R-$tight.\qed
    \end{enumerate}
\end{thm}

\begin{thm}
The following are equivalent for a ring $R$.
    \begin{enumerate}
        \item $R$ is \textit{q.f.d.}.
        \item Every direct sum of minimally tight modules is minimally tight.
        \item Every direct sum of injective $R-$modules is minimally tight.
        \item Every direct sum of indecomposable injectives is minimally tight.\hfill\qed
\end{enumerate}
\end{thm}

Combining Theorem 4.23, Theorem 4.13, Theorem 2.6 in \cite{jain2018survey} and Theorem 3.10 in \cite{kocsan2017weakly}, we get the following result followed by a simpler theorem.

\begin{thm}
    The following are equivalent for a ring $R$.
    \begin{enumerate}
        \item $R$ is \textit{q.f.d.}
        \item Every direct sum of injectives is weakly-injective.
        \item Every direct sum of injectives is tight.
        \item Every direct sum of tight is tight.
        \item Every direct sum of weakly-injectives is tight.
        \item Every direct sum of weakly-injectives is $R-$tight.
        \item Every direct sum of indecomposable injectives is $R-$tight.
        \item Every direct sum of essentially tight is essentially tight.
        \item Every direct sum of injective $R-$modules is essentially tight.
        \item Every direct sum of indecomposable injectives is essentially tight.
        \item Every direct sum of quasi tight modules is quasi tight.
        \item Every direct sum of weakly-injectives is quasi tight.
        \item Every direct sum of weakly-injectives is quasi $R-$tight.
        \item Every direct sum of indecomposable injectives is quasi $R-$tight.
        \item Every direct sum of minimally tight modules is minimally tight.
        \item Every direct sum of weakly-injectives is minimally tight.
        \item Every direct sum of weakly-injectives is minimally $R-$tight.
        \item Every direct sum of indecomposable injectives is minimally $R-$tight.
        \end{enumerate}
\end{thm}

\begin{thm}
   Over a \textit{q.f.d.} ring, the following are equivalent: 
    \begin{enumerate}
        \item $M$ is tight.
        \item $M$ is weakly injective.
        \item $M$ is essentially tight.
        \item $M$ is quasi tight.
        \item $M$ is minimally tight.\hfill \qed
    \end{enumerate}
\end{thm}

\printbibliography

\end{document}